\documentstyle{amsart}
\input{bull-art}
\newtheorem{theo}{Theorem}

\theoremstyle{definition}

\newtheorem{conj}{Conjecture}
\newtheorem{prob}{Problem}

\begin{document}
\def\currentvolume{29}
\def\currentissue{1}
\def\currentyear{1993}
\def\currentmonth{July}
\def\copyrightyear{1993}
\def\currentpages{60-62}

\title{A Counterexample to Borsuk's Conjecture}
\author{Jeff Kahn \and Gil Kalai} 
\address{Department of Mathematics, Rutgers University, 
New Brunswick, New Jersey  08903}
\email{jkahn@@math.rutgers.edu}
\address{Institute of Mathematics,
Hebrew University, Jerusalem  91904,  Israel}
\email{kalai\%humus.huji.ac.il@@relay.cs.net}
\thanks{The first author was
supported in part by BSF, NSF, and AFOSR.
The second author was supported in part by GIF}
\date{June 30, 1992}

\maketitle
\begin{abstract}
Let $f(d)$ be the smallest number so that every set in 
$R^d$ of 
diameter
1 can be partitioned into  $f(d)$ sets of diameter smaller 
than 1.
Borsuk's conjecture was that $f(d)\! =\!d\!+\!1$. 
We prove that $f(d)\! \ge\! (1.2)^{\sqrt d}$
for large~$d$. 
\end {abstract}






\section {Introduction}

Sixty years ago Borsuk \cite {Bo} raised the following 
question.
\begin {prob}[Borsuk]
Is it true that every set of diameter one in $R^d$ 
can be partitioned into $d+1$ closed sets of diameter 
smaller than 
one? 
The conjecture that this is true has come to be called 
Borsuk's 
conjecture.

Let $f(d)$ be the smallest number so that every set in 
$R^d$ of 
diameter
1 can be partitioned into  $f(d)$ sets of diameter smaller 
than 1.
The vertices of the regular simplex in $R^d$ show that 
$f(d)\ge d+1$.
(Another example showing this is, by the 
Borsuk-Ulam theorem, the $d$-dimensional
Euclidean ball.) The assertion of Borsuk's conjecture was 
proved in 
dimensions
2 and 3 and in all dimensions for 
centrally symmetric convex bodies and smooth convex bodies.
See \cite {Gr,BG,CFG} and references cited there. 
Lassak \cite {Las} proved that $f(d) \le 2^{d-1}+1$,
and Schramm \cite {Sc} 
showed that for every $\epsilon$, if $d$ is sufficiently 
large, 
$f(d) \le  (\sqrt{(3/2)}+\epsilon)^d$. 
A different proof of Schramm's bound 
was given by Bourgain and Lindenstrauss \cite {BL}. See 
\cite{Gr,BG,CFG}
for surveys and many references on Borsuk's problem.

Borsuk's conjecture seems to have been believed  
generally, and 
various\break stronger
conjectures have been proposed.  The possibility
of a counterexample based on combinatorial configurations 
was 
suggested by Erd\H{o}s
\cite{A}, Larman \cite{B}, and perhaps others.
In 1965 Danzer \cite{C} showed that the
finite set $K\subseteq{R}^d$ consisting of all
$\{0,1\}$-vectors of an appropriate weight cannot be
covered by $(1.003)^d$ {\em balls} of smaller diameter. 
Larman \cite {La} 
observed that, for sets consisting
of 0-1 vectors with constant weight, Borsuk's conjecture 
reduces to:
\end {prob}
\begin {conj}
Let $K$ be a family of $k$-subsets of $\{1,2,...,n\}$ such 
that every 
two
members of $K$ have $t$ elements in common. Then $K$ can be 
partitioned
into $n$ parts so that in each part every two members  
have $(t+1)$ 
elements
in common.
\end {conj}
Here we prove
\begin {theo} 
For large enough $d$,
$f(d) \ge (1.2)^{\sqrt d}$
by constructing an appropriate family of sets.
\end {theo}

We need the following result of Frankl and Wilson \cite 
{FW}.

\begin {theo} [Frankl and Wilson]
Let $k$ be a prime power and  $n=4k$. Let $K$ be a family
of $n/2$-subsets of $\{1,2,...,n\}$, so that no two
sets in the family have intersection of size $n/4$. Then 
$$
|K| \le 2 \cdot {{n-1} \choose {n/4-1}}. $$
\end {theo}
This settled, in particular, a (much weaker) conjecture of
Larman and Rogers \cite {LR} and implies that, if $g(d)$ 
is the smallest number so that $R^d$ can be colored 
by $g(d)$ colors such that no two points of the same color 
are 
distance one apart, then $g(d)\ge (1.2)^d$.

Let us also mention the following related result
conjectured by Erd\H{o}s and proved by 
Frankl and R\"odl \cite {FR}.

\begin {theo}[Frankl and R\"odl]
Let $n$ be a positive integer divisible by four.
Let $K$ be a family of $n/2$-subsets of $\{1,2,...,n\}$ 
such that no 
two 
sets in the family have intersection of size $n/4$. 
Then $|K|\le (1.99)^{n}$.
\end {theo}

\section{The construction}
\vspace{-1.5pc}
\begin{multline*}
\qquad\mbox{\it However contracted, that definition is the 
result of 
expanded meditation.}\\
\mbox{---Herman Melville, \em Moby Dick}
\end{multline*}


Let $V=\{1,2,...,m\}$, and $m=4k$, and $k$ is a prime power.
Let $W$ be the set of pairs 
of elements in $V$. For every partition $P=\{A,B\}$ of $V$  
let $S(A,B)$ be the sets of all 
pairs which contain one element from $A$ and one element 
from $B$.  
Let $K$ be the family of all 
sets of pairs which correspond to partitions of 
$V$ into two {\it equal} parts, i.e., 
$K=\{S(A,B):|A|=2k\}$. 
Thus, $K$ is a family of ($m^2/4$)-subsets of an 
$m(m-1)/2$-set.
The smallest intersection between $S(A,B)$ and $S(C,D)$ 
occurs
when $|A \cap C|=k$, and by the Frankl-Wilson theorem every 
subfamily
of more than $2 \cdot {{m-1} \choose {m/4-1}}$
sets in $K$ contains two sets which realize
the minimal distance. 
Thus,
$K$ cannot be partitioned into fewer than 
$$
\frac {\frac{1}{2} {{m} \choose {m/2}}} {2 \cdot {{m-1} 
\choose {m/4-1}}}
$$ 
parts so that the 
minimal intersection is not realized in any of the parts.
This expression is greater than $(1.203)^{\sqrt{d}}$ for 
sufficiently
large $d=(^m_2)-1$, and Theorem 1 for general (large) $d$ 
follows
via the prime number theorem. 

\section {Remarks}
1. In view of Theorem 1, the upper bounds on $f(d)$ cited 
earlier 
seem
much more reasonable than formerly.
It would be of considerable 
interest to have a better understanding of the asymptotic
behavior of $\log f(d)$. 
At the moment, we cannot distinguish the asymptotic 
behavior of 
$f(d)$ from that of $g(d)$.
Also of interest would be  
counterexamples in small dimensions. Our construction 
shows that 
Borsuk's conjecture is false for $d=1,325$ and for every 
$d> 2,014$.

2. Larman's conjecture for $t=1$ is open and still quite 
interesting, 
in part because of its similarity to the 
Erd\H{o}s-Faber-Lov\'asz
conjecture. 
See \cite {KS,KK} for some discussion and related results.

3. Intersection properties of edge-sets of graphs were 
first studied 
by S\'{o}s; see \cite {SS} and references quoted therein.

\end{document}